\documentclass[11pt,reqno]{amsart}
\usepackage{mathrsfs}
\usepackage{amssymb}
\usepackage{amsmath}

 \newtheorem{thm}{Theorem}[section]
 
 \newtheorem{lem}[thm]{Lemma}
 
 \newtheorem{exam}[thm]{Example}
 \theoremstyle{definition}
 
 \theoremstyle{remark}
 \newtheorem{rem}[thm]{Remark}
 \numberwithin{equation}{subsection}

\begin{document}

\title{Nonsingular Ricci flow on a noncompact manifold in dimension three}

\author{ Li Ma, Anqiang Zhu}

\address{Li Ma, Department of Mathematical Sciences, Tsinghua University,
 Peking 100084, P. R. China}

\email{lma@math.tsinghua.edu.cn}

\begin{abstract} We
consider the Ricci flow $\frac{\partial}{\partial t}g=-2Ric$ on the
3-dimensional complete noncompact manifold $(M,g(0))$ with
non-negative curvature operator, i.e., $Rm\geq 0, |Rm(p)|\rightarrow
0, ~as ~d(o,p)\rightarrow 0.$ We prove that the Ricci flow on such a
manifold is nonsingular in any finite time.
\end{abstract}



\maketitle

\section{introduction}
The aim of this paper is to get a global existence of Ricci flow
with bounded non-negative curvature operator in three dimensions.
This kind of question was asked by Hamilton \cite{HS}. We remark
that the local existence of the flow was obtained by Shi \cite{S}.
So we only need to show that the curvature is bounded in finite
time. Our research is based on previous important results obtained
by Hamilton and Perelman (\cite{HS}, \cite{KL}, and \cite{MT}),
which will be recalled in next section.

The Ricci flow
$$
\frac{\partial}{\partial t}g_{ij}=-2R_{ij}
$$
was first introduced by Richard Hamilton \cite{H82}. Using it,
Hamilton has obtained some interesting theorem, such as \cite{H82}.
Hamilton's program is to prove Poincar$\acute{e}$ conjecture and
Thurston's geometrization conjecture by Ricci flow. In three
remarkable papers \cite{P1},\cite{P2},\cite{P3}, Perelman
significantly advanced the theory of the Ricci flow. Perelman
introduced canonical neighborhood, and analyze the high curvature
region. Perelman also analyzed one of the special solution to the
Ricci flow, $\kappa$ solution, which is usually the limit solution
of the blow up sequence. Before Perelman, Hamilton\cite{HS} had
defined asymptotic volume and obtained that the asymptotic volume is
constant under Ricci flow. By a induction argument, Perelman
obtained that the asymptotic volume is zero when the solution is
$\kappa$ solution. In order to analyze the high curvature region,
Hamilton obtained some compactness of Ricci flow \cite{HC}. But in
order to apply this compactness, one has to check the non-collapse
and curvature bound assumption. We can use these ideas to prove the
existence of a non-singular Ricci flow on a 3-manifold.

Our main result is the following
\begin{thm}\label{thm4}
Assume that $(M,g(t)), t\in [0,T)$ is a Ricci flow on 3-dimensional
complete noncompact Riemannian manifold. Suppose $Rm(0)\geq 0,
|Rm(p,0)|\rightarrow 0, d(o,p)\rightarrow \infty.$ Then $T=\infty,$
i.e Ricci flow is nonsingular in finite time on such a manifold.
\end{thm}
We now give a interesting example.

\begin{exam}
Consider the paraboloid of revolution
$x_{4}=x_{1}^{2}+x_{2}^{2}+x_{3}^{2},$ $~(x_{1},\cdots ,x_{4})\in
R^{4}.$ We know the curvature $Rm(x)\rightarrow 0, ~x\rightarrow
\infty.$ and $Rm>0$. So Ricci flow on it will not blow up in finite
time. We also know that the asymptotic volume is $0$, so the
paraboloid can't converge uniformly to flat $R^{3}$. We also note
that when the infinite of the manifold is asymptotical to a cone,
the Ricci flow doesn't blow up in finite time. But this time the
asymptotic volume is between $0$ and $\omega,$ where $\omega$ is the
volume of standard ball $B(0,1)\subset R^{3}$.
\end{exam}

We remark that in the radial symmetrical case, a similar result was
obtained in \cite{OW}, where another assumption such as the
asymptotic flatness was used. One may also see our previous work
\cite{M} for more results in this direction.

\section{Some famous results}
In this section, we recall the following results of Hamilton and
Perelman, which will be used in our proof of Theorem \ref{thm4}.

By monotonicity of $\mathscr{W}$ functional and reduced volume,
Perelman proved twice the non-collapse of Ricci flow. Furthermore
Perelman obtained the following convergence theorems about Ricci
flow, which we will use many times (see \cite{P1}, \cite{KL}, and
\cite{MT}).

\begin{thm}\label{thm1}
Fix canonical neighborhood constance $(C,\epsilon),$ and a
non-collapsing constance $r>0,\kappa>0.$ Let
$(\mathscr{M}_{n},G_{n},x_{n})$ be a sequence of based generalized
3-dimensional Ricci flows. We set $t_{n}=t(x_{n})$ and
$Q_{n}=R(x_{n}).$ We denote by $M_{n}$ the time $t_{n}$ time slice
of $\mathscr{M}_{n}.$ We suppose that:
\begin{enumerate}
  \item Each $(\mathscr{M}_{n},G_{n})$ has time interval of
  definition contained in $[0,\infty)$ and has curvature pinched
  toward positive.
  \item Every point $y_{n}\in (\mathscr{M}_{n},G_{n})$ with $t(y_{n})\leq
  t_{n}$ and $R(y_{n})\geq 4R(x_{n})$ has a strong $(C,\epsilon)$
  canonical neighborhood.
  \item $lim_{n\rightarrow \infty} Q_{n}=\infty.$
  \item For each $A<\infty$ the following holds for all n
  sufficiently large. The Ball
  $B(x_{n},t_{n},AQ_{n}^{-\frac{1}{2}})$ has compact closure in
  $M_{n}$ and the flow is $\kappa$ non-collapsed on scales $\leq r$
  at each point of $B(x_{n},t_{n},AQ_{n}^{-\frac{1}{2}}).$
  \item There is $\mu>0$ such that for every $A<\infty$ the
  following holds for all n sufficiently large. For $y_{n}\in
  B(x_{n},t_{n},AQ^{-\frac{1}{2}})$ the maximum flow line through
  $y_{n}$ extends backwards for a time at least $\mu (max(Q_{n},R(y_{n})))^{-1}.$

  Then after passing to a subsequence and shifting the times of each
  the Ricci flow so that $t_{n}=0$ for every n, there is a geometry
  limit $(M_{\infty},g_{\infty},x_{\infty})$ of the sequence of
  based Riemannian manifolds $(M_{n},Q_{n}G_{n}(0),x_{n}).$ The
  limit is a complete 3-dimensional Riemannian manifold of bounded,
  non-negative curvature. Furthermore, for some $t_{0}>0$ which
  depends on the curvature bound for $(M_{\infty},g_{\infty}),$
  there is a geometric limiting Ricci flow defined
  on $(M_{\infty},g_{\infty}(t)), -t_0\leq t\leq
  0,$ with $g_{\infty}(0)=g_{\infty.}$
\end{enumerate}
\end{thm}
\begin{thm}\label{thm2}
Suppose that $\{\mathscr{M}_{n},G_{n},x_{n}\}_{n=1}^{\infty}$ is a
sequence of 3-dimensional Ricci flow satisfying all the hypothesis
of above theorem. Suppose in addition that there is $T_{0}$ with
$0<T_{0}\leq \infty$ such that the following holds. For any
$T<T_{0},$ for each $A<\infty,$ and all n sufficiently large, there
is an embedding $B(x_{n},t_{n},AQ_{n}^{-\frac{1}{2}})\times
(t_{n}-TQ_{n}^{-1},t_{n}]$ into $\mathscr{M}_{n}$ compatible with
time and with the vector field and at every point of the image the
flow is $\kappa$ non-collapsed on scales $\leq r.$ Then, after
shifting the times of the generalized flows so that $t_{n}=0$ for
all n and passing to a subsequence there is a geometric limiting
Ricci flow
$$
(M_{\infty},g_{\infty}(t),x_{\infty}), -T_{0}<t\leq 0,
$$
for the rescaled flows $(\mathscr{M}_{n},Q_{n}G_{n},x_{n}).$ This
limiting flow is complete and of nonnegative curvature. Furthermore,
the curvature is locally bounded n time. If in addition
$T_{0}=\infty,$ then it is a $\kappa$ solution.
\end{thm}

In this paper, we also need the following result of Hamilton on
Ricci flow. By using barrier functions, Hamilton has proved the
asymptotically property is preserved under Ricci flow (see
\cite{HS}).

\begin{thm}\label{thm3}
Assume that we have a solution to the Ricci flow on a complete
noncompact manifold with bounded curvature. If $|Rm(p,0)|\rightarrow
0, d_{g(0)}(o,p)\rightarrow \infty,$ this remains true for $t\geq
0.$
\end{thm}

\begin{rem} One may see \cite{M} for an improvement of this result.
\end{rem}

\section{proof of Theorem \ref{thm4}}
In the following, we will consider the 3-dimensional Ricci flow on
$(M,g(0)),$ with $Rm\geq 0,~|Rm(p)|\rightarrow 0, ~as
~d(o,p)\rightarrow \infty,$ where $o$ is a fixed point. In order to
applied the above convergence theorem, we first show two lemmas. The
method is from Perelman's famous papers (see \cite{P1}, \cite{P2}).
But the condition is different from us here.

\begin{lem}\label{lem1}
For sufficiently small $r>0,$ there is $\kappa>0$ such that the
Ricci flow on this noncompact manifold is $\kappa$ non-collapse on
scales $\leq r.$
\end{lem}

\begin{proof} Fix $(x,t_{0})\in M\times [0,T).$ Since $Rm\geq 0$, Shi has
proved that it's preserved by Ricci flow \cite{S}. By using a well
known result of Gromoll and Meyer \cite{CE}, we have an injectivity
radius estimate
$$
inj(M^{n},g(t))\geq \frac{\pi}{\sqrt{R_{max}(t)}}.
$$
Since
$$
R(x,t)<C_{1}, ~(x,t)\in M\times [0,\frac{1}{2}t_{0}],
$$
by the above injective estimate, we have
$$
Vol B(x,t,r)\geq V^{'}r^{3}, ~(x,t)\in M\times [0,\frac{1}{2}t_{0}].
$$
where $V^{'}$ is a constant. By the inequality of reduced length,
$$
\frac{\partial l_{x}}{\partial \tau}(q,\tau)+\Delta
l_{x}(q,\tau)\leq \frac{(\frac{n}{2})-l_{x}(q,\tau)}{\tau},
$$
we know there is a point $(\tilde{q},\tilde{t}),
~\tilde{t}=\frac{1}{4}t_{0},$ such that
$l_{x}(\tilde{q},\tilde{\tau})\leq \frac{3}{2},$ where
$\tilde{\tau}=t_{0}-\tilde{t}$. By the inequality (see \cite{P1} and
\cite{MT})
$$
|\nabla l_{x}(q,\tau)|^{2}\leq |\nabla
l_{x}(q,\tau)|^{2}+R(q,\tau)\leq \frac{(1+2n)l_{x}(q,\tau)}{\tau},
$$
we have that
$$
l_{x}(q,\tilde{\tau})\leq
(\frac{\sqrt{2n+1}d_{g(t_{0}-\tilde{\tau})}(q,\tilde{q})}{2}+\sqrt{\frac{n}{2}})^{2},
$$
So, for any $A<\infty,$ we have $l_{x}(q,\tilde{\tau})<C(A),$ when
$(q,t_{0}-\tilde{\tau})\in B(\tilde{q},t_{0}-\tilde{\tau},A),$ where
$C(A)$ is a constant depend on $r.$ By Perelman's non-collapse
theorem, we know that if
$$
|Rm(p,t)|<r^{-2}, ~(p,t)\in B(x,t_{0},r)\times [t_{0}-r^{2},t_{0}]
$$
then
$$
Vol B(x,t_{0},r)\geq \kappa r^{n}.
$$
\end{proof}

\begin{lem}\label{lem2}
Fix $0<\epsilon<1.$ Then there is $r>0$ such that for any point
$(x_{0}, t_{0})$ in the flow with $R(x_{0},t_{0})\geq r^{-2}$ the
following hold. $(x_{0},t_{0})$ has a strong canonical
$(C(\epsilon),\epsilon)$ neighborhood.
\end{lem}
\begin{proof} By Hamilton's theorem \ref{thm3}, we know that for any $t^{'}_{n}<T,
~t^{'}_{n}\rightarrow T,$ the curvature is bounded
$$
|Rm(x,t)|<C(t^{'}_{n}), ~(x,t)\in M\times [0,t^{'}_{n}].
$$
Set
$$
A_{n}=\{(x,t)\in M\times [0,t^{'}_{n}]|~(x,t) ~doesn't ~have
~canonical~ neighborhood\}.
$$
Then there is also a up bound to curvature of point in $A_{n},$
$$
R(x,t)<\tilde{C}(t^{'}_{n}), ~(x,t)\in A_{n}.
$$
We can pick point $(x_{n},t_{n})\in A_{n},$ such that
$$
R(x_{n},t_{n})>\frac{1}{2}\tilde{C}(t^{'}_{n}), ~t_{n}\leq t^{'}_{n}
$$ We need to prove
that
$$
\overline{lim}_{n\rightarrow \infty}R(x_{n},t_{n})<\infty.
$$
Suppose not. Then we may assume that
$Q_{n}=R(x_{n},t_{n})\rightarrow \infty.$ From the construction, we
know that $\forall (x,t)\in M\times [0,t_{n}],$ if
$R(x,t)>4R(x_{n},t_{n}),$ $(x,t)$ have a canonical neighborhood.
Since it's a 3-dimensional Ricci flow, the curvature is pinching
toward positive. By Lemma \ref{lem1}, we have the non-collapse
assumption. Since $Q_{n}\rightarrow \infty, ~t_{n}\rightarrow T,$
for any fix $T>0,$ we have $(t_{n}-TQ_{n}^{-1},t_{n}]\subset
[0,t_{n}]$ for sufficiently large n. That is the addition assumption
of Theorem \ref{thm2} is satisfied.

By Theorem \ref{thm2}, we have that
$(M,Q_{n}g(t_{n}+\frac{t}{Q_{n}}),(x_{n},t_{n}))$ converges to a
limit flow $(M_{\infty}, g_{\infty}(t),(x_{\infty},0)),$ which is a
$\kappa$ solution. So for sufficiently large n, $(x_{n},t_{n})$ has
a canonical neighborhood. This contradicts our assumption that none
of the point $(x_{n},t_{n})$ has a canonical neighborhood.
\end{proof}

\emph{ Proof of Theorem \ref{thm4}}: Suppose that the Ricci flow
blows up at time $T.$ By Theorem \ref{thm3}, we know there is a
limit metric $g(T)$ at infinity such that $g(t)|_{M-K}\rightarrow
g(T)|_{M-K},$ where $K$ is a suitable compact set of $M.$ Since the
Ricci flow blows up at the time $t=T,$ we have $sup_{x\in
M}Rm(x,t)\rightarrow \infty, ~t\rightarrow T.$ Otherwise, we can
extend the flow keeping the curvature bounded (That contradicts the
maximum of existence time). So we have that there is a point $p\in
K,$ such that the scalar curvature blows up at $T$, that is
$R(p,t)\rightarrow \infty, ~t\rightarrow T.$ We pick a sequence
$t_{n}\rightarrow T,$ such that $Q_{n}=R(p,t_{n})\rightarrow
\infty.$ By Lemma \ref{lem2}, we have the assumptions of Theorem
\ref{thm2} is satisfied. So
$(M,Q_{n}g(t_{n}+\frac{t}{Q_{n}}),(p,t_{n}))$ converges to a
$\kappa$ solution $(M_{\infty},g_{\infty}(t),(x_{\infty},0)),$

we know that the asymptotic volume of $\kappa$ solution is zero,
that is,
$$lim_{r\rightarrow \infty}\frac{vol
B(x_{\infty},r)}{r^{3}}=0.
$$
Fixing any $\epsilon>0,$ there is a sufficient large $r$ such that
$$
\frac{vol B(x_{\infty},r)}{r^{3}}\leq \frac{\epsilon}{2}.
$$
Since $(M,Q_{n}g(t_{n}+\frac{t}{Q_{n}}),(p,t_{n}))$ converges to
$(M_{\infty},g_{\infty}(t),(x_{\infty},0)),$  for large n, we have
$$
\frac{volB_{Q_{n}g(t_{n})}((p,t_{n}),r)}{r^{3}}\leq \epsilon.
$$
That is
$$
\frac{volB_{g(t_{n})}((p,t_{n}),\frac{r}{Q_{n}^{\frac{1}{2}}})}{(\frac{r}{Q_{n}^{\frac{1}{2}}})^{3}}\leq
\epsilon.
$$
On the other hand, there is a compact region $\Omega\subset M-K,$
such that $g(t)|_{\Omega}\rightarrow g(T)|_{\Omega}.$ Since $R>0,$
$~\frac{d}{dt}\int_{\Omega}d\mu=-\int_{\Omega}Rd\mu\leq 0,$ we have
$$
Vol_{g(t)} \Omega\geq \delta,~t\in [0,T].
$$
Since $\Omega$ is compact, we can find $\tilde{r}>0,$ such that
$$
\Omega \subset B_{g(0)}(p,\tilde{r}).
$$
Since $Ric\geq 0,$ the distance is decreasing. We have $\Omega
\subset B_{g(t)}(p,\tilde{r}), ~t\in [0,T)$. So
$$
Vol B_{g(t_{n})}(p,\tilde{r})\geq \delta.
$$
Since $Q_{n}\rightarrow \infty,$ for sufficiently large $n$,
$\tilde{r}>\frac{r}{Q_{n}}$. We now choose
$\epsilon<\frac{\delta}{\tilde{r}^{3}}.$ Using $Ric\geq 0,$ we have
by the volume comparison theorem that
\begin{eqnarray*}
\frac{\delta}{\tilde{r}^{3}}>\epsilon
>\frac{VolB_{g(t_{n})}((x_{n},t_{n}),\frac{r}{Q_{n}^{\frac{1}{2}}})}{(\frac{r}{Q_{n}^{\frac{1}{2}}})^{3}}
>\frac{Vol B_{g(t_{n})}(x_{n},\tilde{r})}{\tilde{r}^{3}}\geq
\frac{\delta}{\tilde{r}^{3}},
\end{eqnarray*}
which is a contradiction. So the Ricci flow on $(M,g(0))$ is
nonsingular at finite time. That is, $T=\infty.$

\begin{rem}
X. Dai and L. Ma proved that the Ricci flow on the asymptotically
flat manifold can not converge uniformly to flat manifold by ADM
mass (see \cite{M}).
\end{rem}

{\em Acknowledgement:} The research is partially supported by the
National Natural Science Foundation of China 10631020 and SRFDP
20060003002. Part of the work is done when the first name author is
visiting Chern Institute of Mathematics in Nankai University,
Tianjin, which the first named author would like to thank the
director Prof.Y.M.Long for invitation.

\end{document}